\documentclass[12pt]{article}

\usepackage[T2A]{fontenc}
\usepackage[cp1251]{inputenc}
\usepackage[english,russian]{babel}
\usepackage[tbtags]{amsmath}
\usepackage{amsfonts,amssymb}

\usepackage{mathrsfs}

\usepackage{graphicx}

\begin{document}

\begin{center}
M.\,Yu.~Ignatiev

\textbf{On the Weyl--type solutions for differential systems with a singularity. Case of discontinuous potential.}
\end{center}

\section{Introduction}
\label{sec1}

We consider the differential system
\begin{equation}\label{main sys}
  y'-x^{-1}Ay-q(x)y=\rho By, \ x>0,
\end{equation}
with $n\times n$ matrices $A,B, q(x), x\in(0,\infty)$, where $A,B$ are constant, the matrix function $q(\cdot)$ will be referred to as the \textit{potential}.

Differential equations with coefficients having non-integrable singularities at the end or inside the interval often appear in various areas of natural sciences and engineering. As to the case of $n=2$, there exists an extensive literature devoted to different aspects of spectral theory of the radial Dirac operators, see, for instance {\cite{TKBDir}}, {\cite{AHM}}, {\cite{AHM1}}, {\cite{Ser}}, {\cite{GoYu}}. Although this case is very important and has many applications, its investigation do not form a comprehensive picture.
Systems of the form {\eqref{main sys}} with $n>2$ and noncolinear complex eigenvalues of the matrix $B$ show much more complicated behavior and are considerably more difficult for investigation even in the "regular" case $A=0$ {\cite{BCsyst}}.  Some difficulties of principal matter also appear due to presence of the singularity at $x=0$. Whereas the "regular" case $A=0$ has been studied fairly completely to date {\cite{BCsyst}}, {\cite{Zh}}, {\cite{Ysyst}}, for system {\eqref{main sys}} with $A\neq 0$  there are no similar general results.

In this paper, we concentrate on construction and investigation of a distinguished basis of generalized eigenfunctions for {\eqref{main sys}}. We call them  the \textit{Weyl-type solutions}. The Weyl-type solutions play a central role in studying both direct and inverse spectral problems (see, for instance, {\cite{BeDT}},{\cite{Ybook}}). In presence of the singularity at $x=0$ this step encounters some difficulties that do not appear in the "regular" case $A=0$. The approach presented in {\cite{YIP93}} for
the case of closely connected with {\eqref{main sys}} scalar differential operators
\begin{equation}\label{scala}
\ell y=y^{(n)}+\sum\limits_{j=0}^{n-2} \left(\frac{\nu_j}{x^{n-j}}+q_j(x)\right)y^{(j)}
\end{equation}
is based on some special solutions of the equation $\ell y=\lambda y$ that also satisfy certain Volterra integral equations. On this way various results concerning both direct and inverse spectral problems for {\eqref{scala}} were obtained (see, for instance, {\cite{YIP93}}, {\cite{YurMSb}}, {\cite{YurMNinside}} and also {\cite{Fed}} and references therein).  But the approach has some important restriction, namely, it assumes some additional decay condition for the coefficients $q_j(x)$ as $x\to 0$. In this paper, we use another approach that allows us to not impose any additional restrictions of such a type.

In brief outline the approach can be described as follows. We consider some auxiliary systems with respect to the functions with  values  in the exterior algebra $\wedge \mathbb{C}^n$. Our studying of these auxiliary systems centers on two families of  their solutions that also satisfy some asymptotical conditions as $x\to 0$ and $x\to\infty$ respectively, and can be constructed as solutions of certain \textit{Volterra} integral equations. We call these distinguished tensor solutions  the \textit{fundamental tensors}.
The main difference from the above-mentioned method used in {\cite{YIP93}} is that we use the integral equations to construct the fundamental tensors rather than the solutions for the original system. Since each of the fundamental tensors has minimal growth (as $x\to 0$ or $x\to\infty$) among the solutions of the same auxiliary system,
this step does not require any decay of $q(x)$ as $x\to 0$. As a next step, we show that the fundamental tensors are decomposable. Moreover, they can be represented as the wedge products of some solutions of original system {\eqref{main sys}} and these solutions can be shown to be the Weyl-type solutions of {\eqref{main sys}}. Such strategy allows us to construct the Weyl-type solutions via purely algebraic procedure and investigate their properties.

This approach was first presented in {\cite{BeDT}}
for the higher-order differential operators with {\it regular} (from the Schwarz space) coefficients on the whole line and was used mostly for evaluating a behavior of the Weyl--type solutions and of the scattering data for $\rho\to 0$. In {\cite{IgnSingRM16}} the approach was modified and adapted to the case of system {\eqref{main sys}}. The developed method was used for investigation of asymptotical behavior of the Weyl--type solutions as $\rho\to\infty$ and $\rho\to 0$, moreover, some uniqueness result for the inverse scattering problem for {\eqref{main sys}} was also obtained.

In the work {\cite{IgnSingRM16}} the potential $q(\cdot)$ was assumed to be integrable and absolutely continuous on the semi-axis. In this paper, we consider the case when the function $q(\cdot)$ is allowed to be non-differential and even discontinuous. Differential systems with discontinuous coefficients are investigated intensively during last two decades (see, for instance, {\cite{SaShD}}, {\cite{HryZSint}}), one should also mention an important role played by such systems in theory of ordinary differential operators with distribution coefficients (see, for instance, {\cite{MSha}}, {\cite{KMSha}}, {\cite{SaShFSS}}). Investigation of such systems required further nontrivial development of the spectral theory methods (see {\cite{Rykh}}, {\cite{SaCZ}}, {\cite{SaShFSS}}).

The results obtained in the present paper include, in particular, the asymptotics of the Weyl--type solutions as $\rho\to\infty$ with $o(1)$ remainder's  estimates. Results of such a type play an auxiliary role, more important results are formulated in terms of belonging of the Weyl--type solutions to certain classes of integrable functions of spectral parameter $\rho$. We also establish the continuity of  the Weyl--type solutions with respect to the potential $q(\cdot)$.

\section{Solutions of the unperturbed system}
\label{sec2}

In this section we briefly discuss the unperturbed system:
\begin{equation}\label{unp sys}
  y'-x^{-1}Ay=\rho By
\end{equation}
and introduce some fundamental systems of its solutions.

Here and throughout the paper we assume the following

\medskip
\textbf{Condition A.} $A$ is off-diagonal.
The eigenvalues $\{\mu_j\}_{j=1}^n$ of the matrix $A$ are distinct and such that $\mu_j-\mu_k \notin \mathbb{Z}$ for $j\neq k$, moreover, $\mbox{Re}\mu_1<\mbox{Re}\mu_2<\dots<\mbox{Re}\mu_n$ and $\mbox{Re}\mu_k\neq 0$, $k=\overline{1,n}$.

\medskip
\textbf{Condition B.}
$B=diag(b_1,\dots,b_n)$, where the entries $b_1,\dots,b_n$ are nonzero distinct complex numbers such that any three of them are noncolinear and, moreover, $\sum\limits_{j=1}^n b_j=0$.

We start with considering {\eqref{unp sys}} for $\rho=1$:
\begin{equation}\label{un sys 1}
  y'-x^{-1}Ay=By
\end{equation}
but for {\it complex} values of $x$.

Let $\Sigma$ be the following union of lines through the origin in $\mathbb{C}$:
$$
\Sigma=\bigcup\limits_{(k,j): j\neq k}\left\{x:\mbox{Re}(xb_j)=\mbox{Re}(x b_k)\right\}.
$$

Consider some (arbitrary) open sector $\mathcal S\subset\mathbb{C}\setminus\Sigma$ with the vertex at $x=0$. It is well-known that there exists the ordering $R_1,\dots, R_n$ of the numbers $b_1,\dots,b_n$ such that $\mbox{Re}(R_1 x)<\mbox{Re}(R_2 x)<\dots<\mbox{Re}(R_n x)$ for any $x\in \mathcal S$. For $x\in\overline{\mathcal S}\setminus\{0\}$ there exist the following fundamental matrices for system {\eqref{un sys 1}}(see {\cite{CoLe}}, {\cite{Syb}}, {\cite{Fed90}}):
\begin{itemize}
\item $c(x)=(c_1(x),\dots,c_n(x))$, where
$$
c_k(x)=x^{\mu_k}\hat c_k(x),
$$
$\det c(x)\equiv 1$ and all $\hat c_k(\cdot)$ are entire functions, $\hat c_k(0)=\mathfrak{h}_k$, $\mathfrak{h}_k$ is an eigenvector of the matrix $A$ corresponding to the eigenvalue $\mu_k$;
\item $e(x)=(e_1(x),\dots,e_n(x))$, where
$$
e_k(x)=\exp\left(x R_k\right)\left(\mathfrak{f}_k+O\left(x^{-1}\right)\right), \ |x|\geq 1, x\in\overline{\mathcal S},
$$
$(\mathfrak{f}_1,\dots,\mathfrak{f}_n)=\mathfrak f$ is a permutation matrix such that $(R_1,\dots,R_n)=(b_1,\dots,b_n)\mathfrak f$.
\end{itemize}

\medskip
\textbf{ Condition $I(\mathcal S)$.} For all $k=\overline{2,n}$ the Wronskian determinants $$\Delta^0_k:=\det(e_1(x),\dots,e_{k-1}(x),c_k(x),\dots,c_n(x))$$ are not equal to 0.

\medskip
Under Condition $I(\mathcal S)$ for $x\in\overline{\mathcal S}\setminus\{0\}$ there exists (and is unique) the fundamental matrix $\psi_0(x)=(\psi_{01}(x), \dots, \psi_{0n}(x))$ such that
\begin{equation}\label{as psi0}
  \psi_{0k}(xt)=\exp(xt R_k)(\mathfrak{f}_k+o(1)), \ t\to +\infty, x\in\mathcal{S} , \ \psi_{0k}(x)=O(x^{\mu_k}), x\to 0.
\end{equation}
Since $\psi_{0k}(x)$ admits the representation:
$$
\psi_{0k}(x) = \sum\limits_{j=k}^n l_{jk} c_{j}(x) = \sum\limits_{j=k}^n l_{jk} x^{\mu_j} \hat c_{j}(x)
$$
with constant coefficients $\{l_{jk}\}$ we notice that the function $x^{-\mu_k}\psi_{0k}(x)$ admits a continuous extension onto $\overline{\mathcal S}$.

\medskip
Now we return to unperturbed system {\eqref{unp sys}} with arbitrary $\rho\in \overline{\mathcal S}\setminus\{0\}$ and real positive $x$. Notice that if some function $y=y(x)$ satisfies {\eqref{un sys 1}} then $Y(x,\rho):=y(\rho x)$ satisfies {\eqref{unp sys}}. Taking this into account, we define the matrix solutions $C(x,\rho)$, $E(x,\rho)$, $\Psi_0(x,\rho)$ of {\eqref{unp sys}} as follows: $C(x,\rho):=c(\rho x)$, $\Psi_0(x,\rho):=\psi_0(\rho x)$, $E(x,\rho):=e(\rho x)$. The $k$-th column $\Psi_{0k}(x,\rho)$ of the matrix $\Psi_0(x,\rho)$ will be referred to as \textit{$k$-th Weyl type solution} of system {\eqref{unp sys}}. For $\rho\in\mathcal S$ it satisfies the asymptotic conditions:
\begin{equation}\label{un We as}
  \Psi_{0k}(x,\rho) = \exp(\rho x R_k)\left(\mathfrak f_k + o(1)\right), x\to+\infty, \
  \Psi_{0k}(x,\rho) = O\left(x^{\mu_k}\right), x\to 0.
\end{equation}

In the sequel we shall also use the following estimate that can be obtained directly from the definition of the matrix $\Psi_0(x,\rho)$ and {\eqref{as psi0}}:
\begin{equation}\label{est psi0}
  \|\Psi_{0k}(x,\rho)\|\leq M |W_k(\rho x)|.
\end{equation}
Here and below $W_k(\cdot)$, $k=\overline{1,n}$ denotes the weight function defined as follows:
$$
W_k(\xi):=\left\{ \begin{array}{l}
W_0\left(\xi^{\mu_k}\right)\exp(R_k \xi), \ |\xi|\leq 1 \\
\exp(R_k \xi), \ |\xi|> 1,
\end{array}
\right.
$$
$$
W_0(\xi)=(1-|\xi|)\xi+|\xi|^2, \ |\xi|\leq 1, \ W_0(\xi):=(W_0\left(\xi^{-1}\right))^{-1}, \ |\xi|>1,
$$
$M$ is some absolute constant (we agree to use the same symbol $M$ to denote possibly different absolute constants).

From the construction and  the properties of the matrices $\psi_0(x)$ and $\Psi_0(x,\rho)$ one can deduce, moreover, that the functions
$$
\tilde\Psi_{0k}(x,\rho):= (W_k(\rho x))^{-1}\Psi_{0k}(x,\rho)
$$
admit continuous extensions onto the set $\left\{(x,\rho):x\in[0,\infty), \rho\in\overline{\mathcal S}\right\}$.

\section{Fundamental tensors}
\label{sec3}

Proceeding the way described above we consider for $m=\overline{1,n}$ the following auxiliary systems:
\begin{equation}\label{tens sys}
  Y'=(U(x,\rho))^{(m)}Y,
\end{equation}
where $Y$ is a function with values in the exterior product $\wedge^m \mathbb{C}^n$. Here and below
$$
U(x,\rho):=x^{-1}A+q(x)+\rho B
$$
and for $n\times n$ matrix $V$ the symbol  $V^{(m)}$ denotes the operator acting in $\wedge^m \mathbb{C}^n$ so that for any vectors $v_1, \dots, v_m$ the following identity holds:
$$
V^{(m)}(v_1\wedge v_2 \wedge\dots\wedge v_m)=\sum\limits_{j=1}^m v_1\wedge v_2 \wedge\dots\wedge v_{j-1}\wedge V v_j\wedge v_{j+1}\wedge\dots\wedge v_m.
$$

The system $Y'=U^{(m)}Y$ is just the system which is satisfied by $m$-vector $Y=y_1\wedge\dots\wedge y_m$ constructed by solutions $y_1, \dots, y_m$ of the system $y'=Uy$. Systems {\eqref{tens sys}} can be written, of cause, in the similar form as the original system. The idea is to restrict our analysis with only some special particular solutions of {\eqref{tens sys}}, namely the solutions which have a minimal growth (among the solutions of the same auxiliary system) as $x\to 0$ or as $x\to\infty$. The key technical observation consists with the fact that such solutions of {\eqref{tens sys}} can be constructed as solutions of certain \textit{Volterra} integral equations.

In order to proceed a way described we have to introduce some notations.
\begin{itemize}
  \item $\mathcal{A}_{m}$ denotes the set of all ordered multi-indices $\alpha=(\alpha_1, \dots, \alpha_m)$, $\alpha_1<\alpha_2<\dots<\alpha_m$, $\alpha_j\in\{1,2,\dots,n\}$.
  \item For a set of vectors $u_1, \dots, u_n$
 from $\mathbb{C}^n$ and a multi-index $\alpha\in \mathcal{A}_m$ we define
 $$u_\alpha:=u_{\alpha_1}\wedge\dots\wedge u_{\alpha_m}.$$
  \item Let $\{a_1,\dots,a_n\}$ be a numerical sequence. For $\alpha\in\mathcal{A}_m$ we define $$a_\alpha:=\sum\limits_{j\in\alpha} a_j,$$ for $k\in\overline{1,n}$ we denote $$\overrightarrow{a}_k:=\sum\limits_{j=1}^k a_j,\ \overleftarrow{a}_k:=\sum\limits_{j=k}^n a_j.$$
  \item For a multi-index $\alpha$ the symbol $\alpha'$ denotes the ordered multi-index that complements $\alpha$ to $(1,2,\dots,n)$.
  \item For $h\in \wedge^n \mathbb{C}^n$ we define $|h|$ as a constant in the following representation:
$$
h=|h|\mathfrak{e}_1\wedge\mathfrak{e}_2\wedge\dots\wedge\mathfrak{e}_n,
$$
where (and everywhere below) $\mathfrak{e}_1,\mathfrak{e}_2,\dots,\mathfrak{e}_n$ is the standard basis in $\mathbb C^n$.
\item The space $\wedge^m \mathbb{C}^n$ will be considered with $l_1$ norm:
$$
\left\|\sum\limits_{\alpha\in\mathcal{A}_m}h_\alpha \mathfrak{e}_\alpha\right\|:=\sum\limits_{\alpha\in\mathcal{A}_m}|h_\alpha| .
$$
\end{itemize}

As above we fix some arbitrary open sector $\mathcal S \subset\mathbb C\setminus\Sigma$.

Consider the following Volterra integral equations:
\begin{equation}\label{eq Tk}
  Y(x)=T_k^0(x,\rho)+\int\limits_0^x G_{n-k+1}(x,t,\rho)\left(q^{(n-k+1)}(t)Y(t)\right)dt,
\end{equation}
\begin{equation}\label{eq Fk}
  Y(x)=F_k^0(x,\rho)-\int\limits_x^\infty G_{k}(x,t,\rho)\left(q^{(k)}(t)Y(t)\right)dt,
\end{equation}
where
$$T_k^0(x,\rho):=C_k(x,\rho)\wedge\dots\wedge C_n(x,\rho),$$
$$F_k^0(x,\rho):=E_1(x,\rho)\wedge\dots\wedge E_k(x,\rho)=\Psi_{01}(x,\rho)\wedge\dots\wedge \Psi_{0k}(x,\rho)$$
and $G_m(x,t,\rho)$ is an operator acting in $\wedge^m \mathbb{C}^n$ as follows:
$$
G_m(x,t,\rho)f=\sum\limits_{\alpha\in\mathcal{A}_m}\sigma_\alpha\left|f\wedge C_{\alpha'}(t,\rho)\right|C_{\alpha}(x,\rho)=
$$$$
\sum\limits_{\alpha\in\mathcal{A}_m}\chi_\alpha\left|f\wedge\Psi_{0\alpha'}(t,\rho)\right|\Psi_{0\alpha}(x,\rho)=
\sum\limits_{\alpha\in\mathcal{A}_m}\chi_\alpha\left|f\wedge E_{\alpha'}(t,\rho)\right|E_{\alpha}(x,\rho).
$$
Here  $\sigma_\alpha=\left|\mathfrak{h}_\alpha\wedge\mathfrak{h}_{\alpha'}\right|$, $\chi_\alpha=\left|\mathfrak{f}_\alpha\wedge\mathfrak{f}_{\alpha'}\right|$.

Repeating the arguments used in proof of {\cite{IgnSingRM16}}, Theorem 3.1, we obtain the following result.

\medskip
\textbf{Proposition 1.} \textit{
Suppose that $q\in L_1(0,\infty)$. Then
for any $\rho\in\overline{\mathcal S}\setminus\{0\}$ equations {\eqref{eq Tk}} and {\eqref{eq Fk}} have unique solutions $T_k(\cdot,\rho)$ and $F_k(\cdot,\rho)$ respectively such that
$$
\|T_k(x,\rho)\|\leq M \left\{ \begin{array}{l}
\left|(\rho x)^{\overleftarrow{\mu}_k}\right|, \ |\rho x|\leq 1 \\
\left|\exp(\rho x \overleftarrow{R}_k)\right|, \ |\rho x|>1,
\end{array}\right.$$$$
\|F_k(x,\rho)\|\leq M \left\{ \begin{array}{l}
\left|(\rho x)^{\overrightarrow{\mu}_k}\right|, \ |\rho x|\leq 1 \\
\left|\exp(\rho x \overrightarrow{R}_k)\right|. \ |\rho x|>1.
\end{array}\right.
$$
The following asymptotics hold:
$$
F_k(x,\rho)=\exp(\rho x \overrightarrow{R}_k)\left(\mathfrak{f}_1\wedge\dots\wedge\mathfrak{f}_k+o(1)\right), \ x\to\infty,
$$
$$
T_k(x,\rho)=(\rho x)^{\overleftarrow{\mu}_k}\left(\mathfrak{h}_k\wedge\dots\wedge\mathfrak{h}_n+o(1)\right), \ x\to 0.
$$
}

In what follows we consider the fundamental  tensors as functions of \textit{three} arguments $q,x,\rho$, where the potential $q$ will be assumed to belong to
$\mathcal X_p$, where $\mathcal X_p$ denotes the Banach space consisting of all off--diagonal matrix functions with entries from $X_p:=L_1(0,\infty)\cap L_p(0,\infty)$.

For a closed unbounded subset $L$ of the closed sector $\overline{\mathcal S}$ we denote by $C_0(L)$ the Banach space of continuous  on $L$ and vanishing at infinity functions with the standard $\sup$ norm.

\medskip
\textbf{Theorem 1.}\textit{
Suppose $p>2$. Then:}

\textit{
I. The function $T_k(q,x,\rho)$, $q\in \mathcal X_p, x\in(0,\infty), \rho\in \overline{\mathcal S}\setminus\{0\}$ admits the representation:
$$
T_k(q,x,\rho)=T^0_k(x,\rho)+\overleftarrow W^k(\rho x)\hat T_k(q,x,\rho),
$$
where $\hat T_k\in C\left(\mathcal X_p, BC\left([0,\infty), C_0\left(\overline{\mathcal S}\right)\right)\right)$. Moreover, for any ray of the form
  $\{\rho=zt, t\in[0,\infty)\}$ with $z\in\overline{\mathcal{S}}\setminus\{0\}$  the restriction
$\left.\hat T_k\right|_l$ belongs to the space  $C\left(\mathcal{X}_p, BC\left([0,\infty),\mathcal{H}(l)\right)\right)$, where $\mathcal{H}(l):= C_0(l)\cap L_2(l)$.}

\textit{
II. The function $F_k(q,x,\rho)$, $q\in \mathcal X_p, x\in(0,\infty),  \rho\in \overline{\mathcal S}\setminus\{0\}$ admits the representation:
$$
F_k(q,x,\rho)=F^0_k(x,\rho)+\overrightarrow W^k(\rho x)\hat F_k(q,x,\rho),
$$
where $\hat F_k\in C\left(\mathcal X_p, BC\left([0,\infty), C_0\left(\overline{\mathcal S}\right)\right)\right)$. Moreover, for any ray
  $\{\rho=zt, t\in[0,\infty)\}$ with $z\in\overline{\mathcal{S}}\setminus\{0\}$  the restriction
$\left.\hat F_k\right|_l$ belongs to the space $C\left(\mathcal{X}_p, BC\left([0,\infty),\mathcal{H}(l)\right)\right)$.
}

\medskip
\textbf{Proof.}
I. Substitution
 $$
T_k(q,x,\rho)=T^0_k(x,\rho)+\overleftarrow W^k(\rho x)\hat T_k(q,x,\rho)
$$
converts {\eqref{eq Tk}} into the form
$$\hat T_k(q,\cdot,\cdot)=\mathcal K(q)\hat T_k(q,\cdot,\cdot)+\hat T^1_k(q,\cdot,\cdot),$$
where
$$\hat T^1_k(q,x,\rho)=\left(\overleftarrow W^k(\rho x)\right)^{-1}
\int\limits_0^x G_{n-k+1}(x,t,\rho)\left(q^{(n-k+1)}(t)T^0_k(t,\rho)\right)dt,$$
$\mathcal K(q)$ -- linear operator mapping a function $f=f(x,\rho)$ into the function:
\begin{equation}\label{oper K}
  (\mathcal K(q)f)(x,\rho)=\int\limits_0^x \mathcal G_{n-k+1}(x,t,\rho) \left(q^{(n-k+1)}(t)f(t,\rho)\right)\,dt,
\end{equation}
\begin{equation}\label{mathcal G}
  \mathcal{G}_{n-k+1}(x,t,\rho):=\frac{\overleftarrow{W}^k(\rho t)}{\overleftarrow{W}^k(\rho x)}G_{n-k+1}(x,t,\rho).
\end{equation}

From the results of {\cite{IgnIntTr}} it follows that $\hat T^1_k\in C\left(\mathcal X_p, BC\left([0,\infty), C_0\left(\overline{\mathcal S}\right)\right)\right)$ and, moreover, for any ray
  $\{\rho=zt, t\in[0,\infty)\}$ with $z\in\overline{\mathcal{S}}\setminus\{0\}$  the restriction
$\left.\hat T^1_k\right|_l$ belongs to the space $C\left(\mathcal{X}_p, BC\left([0,\infty),\mathcal{H}(l)\right)\right)$.

I. 1) Let us show that the operator $\mathcal K(q)$ defined by {\eqref{oper K}} acts and is continuous in $BC\left([0,\infty), C_0\left({\overline{\mathcal{S}}}\right)\right)$.

We notice first that the function $\mathcal G_{n-k+1}(x,x\tau,\rho)$ is continuous in $(x,\tau,\rho)\in[0,\infty)\times[0,1]\times\overline{\mathcal S}$ and therefore for any
$f(\cdot,\cdot)\in BC\left([0,\infty), C_0\left({\overline{\mathcal{S}}}\right)\right)$ the function
$$
 (\mathcal K(q)f)(x,\rho)=\int\limits_0^1 \mathcal G_{n-k+1}(x,x\tau,\rho) \left(q^{(n-k+1)}(x\tau)f(x\tau,\rho)\right)\,xd\tau
$$
is continuous in $(x,\rho)\in [0,\infty)\times{\overline{\mathcal{S}}}$. Furthermore, from {\eqref{est psi0}} it follows that the function $\mathcal G_{n-k+1}(x,x\tau,\rho)$ is bounded. This yields the estimate:
\begin{equation}\label{est K}
  \left\|(\mathcal K(q)f)(x,\rho)\right\|\leq M \int\limits_0^x \|q(t)\| \|f(t,\rho)\|\,dt.
\end{equation}
From {\eqref{est K}} one can deduce, in particular:
\begin{equation}\label{est K C0}
  \left\|(\mathcal K(q)f)(x,\rho)\right\|\leq M \sup\limits_{t\in[0,x]}\|f(t,\rho)\|\int\limits_0^x \|q(t)\|\,dt.
\end{equation}
For any $f(\cdot,\cdot)\in BC\left([0,\infty), C_0\left({\overline{\mathcal{S}}}\right)\right)$ we have $\lim\limits_{\rho\to\infty,\rho\in{\overline{\mathcal{S}}}}\sup\limits_{t\in[0,x]}\|f(t,\rho)\|=0$ and from {\eqref{est K C0}} it follows that
$(\mathcal K(q)f)(x,\cdot)\in C_0\left({\overline{\mathcal{S}}}\right)$ for any fixed $x$, moreover,
$\sup\limits_{x\in[0,T]}\|(\mathcal K(q)f)(x,\rho)\|\to 0$ as $\rho\to\infty$, $\rho\in{\overline{\mathcal{S}}}$ for any finite $T>0$.
Since $(\mathcal K(q)f)(\cdot,\cdot)\in C\left(([0,\infty)\times\overline{\mathcal S}\right)$ this means that
$(\mathcal K(q)f)(\cdot,\cdot)\in C([0,\infty), C_0({\overline{\mathcal{S}}}))$. Finally, estimate {\eqref{est K C0}} implies the estimate:
\begin{equation}\label{est K C0 norm}
  \left\|(\mathcal K(q)f)(x,\rho)\right\|\leq M \|f\|_{BC([0,\infty), C_0({\overline{\mathcal{S}}}))} \cdot \|q\|_{L_1(0,\infty)}.
\end{equation}
From {\eqref{est K C0 norm}} we conclude that $(\mathcal K(q)f)(\cdot,\cdot)\in BC([0,\infty), C_0({\overline{\mathcal{S}}}))$ and  $\mathcal K(q)\in \mathcal L\left(BC([0,\infty), C_0({\overline{\mathcal{S}}}))\right)$. Moreover, since the mapping $q\to\mathcal K(q)$ is linear, from {\eqref{est K C0 norm}} it follows that the operator $\mathcal K(q)\in \mathcal L\left(BC([0,\infty), C_0({\overline{\mathcal{S}}}))\right)$ depends continuously on $q\in\mathcal X_p$.

I. 2) From {\eqref{est K}}
it follows the estimate for the iterated operators $\mathcal K^r(q)$:
\begin{equation}\label{est Kr}
  \left\|(\mathcal K^r(q)f)(x,\rho)\right\|\leq \frac{M^r}{r!} \left(\int\limits_0^x \|q(t)\|\,dt\right)^r
  \sup\limits_{\tau\in[0,x]}\|f(\tau,\rho)\|.
\end{equation}
This yields:
\begin{equation}\label{est Kr norm}
  \left\|\mathcal K^r(q)\right\|\leq \frac{M^r}{r!} \|q\|^r_{L_{1}(0,\infty)},
\end{equation}
where norm in the left hand side assumes the $\mathcal L\left(BC([0,\infty), C_0({\overline{\mathcal{S}}}))\right)$ norm.
Therefore, the operator $Id-\mathcal K(q)$ is invertible in $BC([0,\infty), C_0({\overline{\mathcal{S}}}))$ for any $q\in\mathcal X_p$, moreover, the operator $\left(Id-\mathcal K(q)\right)^{-1}$ depends continuously on $q$. Since
$$\hat T_k(q)=(Id-\mathcal K(q))^{-1}\hat T^1_k(q),$$
this means that $\hat T_k\in C\left(\mathcal X_p, BC\left([0,\infty), C_0\left({\overline{\mathcal{S}}}\right)\right)\right)$.

I. 3) Consider an arbitrary ray $l=\{\rho=zt, t\in[0,\infty)\}$ with $z\in\overline{\mathcal{S}}\setminus\{0\}$. For arbitraty $f\in BC([0,\infty), \mathcal H(l))$ the arguments from I. 1) are still valid and yield, in particular, that $(\mathcal K(q)f)(x,\rho)$ is continuous in $(x,\rho)\in[0,\infty)\times l$. Moreover, for $\rho\in l$ inequalities {\eqref{est K}}, {\eqref{est K C0}} are true and therefore $(\mathcal K(q)f)(\cdot,\cdot)\in BC([0,\infty), C_0(l))$.

In what follows we denote by $l^+(R)$ the open ray $l\cap\{\rho:|\rho|>R\}$. Using the estimates:
\begin{equation}\label{est K L2l}
  \left\|(\mathcal K(q)f)(x,\cdot)\right\|_{L_2(l)}
  \leq M \int\limits_0^x \|q(t)\|\|f(t,\cdot)\|_{L_2(l)} \, dt,
\end{equation}
$$
  \left\|(\mathcal K(q)f)(x,\cdot)\right\|_{L_2(l^+(R))}
  \leq M \int\limits_0^x \|q(t)\|\|f(t,\cdot)\|_{L_2(l^+(R))} \, dt,
$$
with arbitrary $R>0$ and Lemma 3.1 from {\cite{IgnIntTr}} we deduce that $(\mathcal K(q)f)(x,\cdot)\in L_2(l)$ for any $x\in[0,\infty)$ and, moreover, $(\mathcal K(q)f)(\cdot,\cdot)\in BC\left([0,\infty), L_2(l)\right)$. From {\eqref{est K C0}}, {\eqref{est K L2l}} it follows also:
$$
 \left\|(\mathcal K(q)f)(x,\rho)\right\|\leq M \|f\|_{BC([0,\infty), C_0(l))} \cdot \|q\|_{L_1(0,\infty)}, $$$$
   \left\|(\mathcal K(q)f)(x,\cdot)\right\|_{L_2(l)}
  \leq M \|f\|_{BC([0,\infty), L_2(l))} \cdot \|q\|_{L_1(0,\infty)},
$$
that allows us to conclude that $\mathcal K(q)\in \mathcal L\left(BC([0,\infty), \mathcal H(l))\right)$ and depends continuously on $q\in\mathcal X_p$.

I. 4) From {\eqref{est K L2l}} we deduce the estimate:
$$
\left\|(\mathcal K^r(q)f)(x,\cdot)\right\|_{L_2(l)}
  \leq \frac{M^r}{r!} \left(\int\limits_0^x \|q(t)\|\, dt\right)^r \sup\limits_{\tau\in[0,x]}\|f(\tau,\cdot)\|_{L_2(l)},
$$
that, together with {\eqref{est Kr}} imply the estimate similar to
{\eqref{est Kr norm}} but with the $\mathcal L\left(BC([0,\infty), \mathcal H(l))\right)$ norm in the left hand side. Proceeding as in I. 2), we conclude that $\hat T_k\in C\left(\mathcal X_p, BC\left([0,\infty), \mathcal H(l)\right)\right)$.

II. 1) Proceeding as in I. 1) we consider $\hat F_k(q,\cdot,\cdot)$ as a solution of the equation $\hat F_k(q)=\mathcal K^+(q)\hat F_k(q)+\hat F^1_k(q)$, where
\begin{equation}\label{oper K plus}
  (\mathcal K^+(q)f)(x,\rho)=-\int\limits_x^\infty \mathcal G_{k}(x,t,\rho) \left(q^{(k)}(t)f(t,\rho)\right)\,dt,
\end{equation}
while the free term
$\hat F^1_k(q)$ can be evaluated by using the results of {\cite{IgnIntTr}}. That yields:
$$\hat F^1_k\in C\left(\mathcal X_p, BC\left([0,\infty), C_0\left({\overline{\mathcal{S}}}\right)\right)\right)$$
and for any ray
  $\{\rho=zt, t\in[0,\infty)\}$ with $z\in\overline{\mathcal{S}}\setminus\{0\}$  the restriction
$\left.\hat F^1_k\right|_l$ belongs to \\ $C\left(\mathcal{X}_p, BC\left([0,\infty),\mathcal{H}(l)\right)\right)$.

Since the function $\mathcal G_{k}(x,t,\rho)$ is continuous in $(x,t,\rho)$: $0\leq x<t<\infty$, $\rho\in{\overline{\mathcal{S}}}$, for any  $f(\cdot,\cdot)\in BC\left([0,\infty), C_0\left({\overline{\mathcal{S}}}\right)\right)$ the function $(\mathcal K^+(q)f)(x,\rho)$ is continuous in $(x,\rho)\in[0,\infty)\times{\overline{\mathcal{S}}}$.
From the boundedness of the function $\mathcal G_{k}(x,t,\rho)$ it follows the estimate :
\begin{equation}\label{est K plus}
  \left\|(\mathcal K^+(q)f)(x,\rho)\right\|\leq M \int\limits_x^\infty \|q(t)\| \|f(t,\rho)\|\,dt.
\end{equation}
Let us show that:
\begin{equation}\label{for K plus to C0}
  \lim\limits_{\rho\to\infty}\sup\limits_{x\in[0,\infty)}\left\|(\mathcal K^+(q)f)(x,\rho)\right\|=0.
\end{equation}
Fix arbitrary $\varepsilon>0$ and choose $T_*=T_*(\varepsilon)$ such that
$$
M\|f\|\int\limits_{T_*}^\infty \|q(t)\| \,dt < \frac{\varepsilon}{2},
$$
where $\|f\|=\|f\|_{BC([0,\infty), C_0({\overline{\mathcal{S}}}))}$ and $M$ is the same as in {\eqref{est K plus}}. For chosen (finite) $T_*$ one can find $R$ such that
$$
M\|f(x,\rho)\|\int\limits_0^{T_*} \|q(t)\| \,dt < \frac{\varepsilon}{2}
$$
for all $x\in[0,T_*]$, $\rho\in{\overline{\mathcal{S}}}: |\rho|>R$. Then for $x\in[0,T_*]$, $\rho\in{\overline{\mathcal{S}}}: |\rho|>R$ we have $\left\|(\mathcal K^+(q)f)(x,\rho)\right\|<\varepsilon$, while for  $x\in[T_*,\infty]$, $\rho\in{\overline{\mathcal{S}}}$ we have $\left\|(\mathcal K^+(q)f)(x,\rho)\right\|<\varepsilon/2$. Thus, {\eqref{for K plus to C0}} is established. Since the function $(\mathcal K^+(q)f)(x,\rho)$ is continuous in $x\in[0,\infty)$, $\rho\in{\overline{\mathcal{S}}}$ from {\eqref{for K plus to C0}} it follows that $(\mathcal K^+(q)f)(\cdot,\cdot)\in BC\left([0,\infty), C_0\left({\overline{\mathcal{S}}}\right)\right)$.

Finally, {\eqref{est K plus}} implies the estimate
$$
\left\|(\mathcal K^+(q)f)(x,\rho)\right\|\leq M  \|q\|_{L_1(0,\infty)} \|f\|_{BC\left([0,\infty), C_0\left({\overline{\mathcal{S}}}\right)\right)},
$$
from which we deduce that the operator $\mathcal K^+(q)\in \mathcal L\left(BC\left([0,\infty), C_0\left({\overline{\mathcal{S}}}\right)\right)\right)$ and continuous with respect to $q\in\mathcal X_p$.

II. 2) From {\eqref{est K plus}}
we deduce the estimate:
\begin{equation}\label{est Kr plus}
  \left\|((\mathcal K^+(q))^r f)(x,\rho)\right\|\leq \frac{M^r}{r!} \left(\int\limits_x^\infty \|q(t)\|\,dt\right)^r
  \sup\limits_{\tau\in[x,\infty)}\|f(\tau,\rho)\|,
\end{equation}
which yields:
\begin{equation}\label{est Kr plus norm}
  \left\|(\mathcal K^+(q))^r\right\|\leq \frac{M^r}{r!} \|q\|^r_{L_{1}(0,\infty)},
\end{equation}
where the norm in the left hand side assumes the $\mathcal L\left(BC([0,\infty), C_0({\overline{\mathcal{S}}}))\right)$ norm.
Thus, one can conclude that the operator $Id-\mathcal K^+(q)$ is invertible in $BC([0,\infty), C_0({\overline{\mathcal{S}}}))$ for any $q\in\mathcal X_p$, moreover, the operator $\left(Id-\mathcal K^+(q)\right)^{-1}$ is continuous with respect to $q$. Therefore, we have proved that $\hat F_k\in C\left(\mathcal X_p, BC\left([0,\infty), C_0\left({\overline{\mathcal{S}}}\right)\right)\right)$.

II. 3) Consider an arbitrary ray $l=\{\rho=zt, t\in[0,\infty)\}$ with $z\in\overline{\mathcal{S}}\setminus\{0\}$. For arbitrary $f\in BC([0,\infty), \mathcal H(l))$  the arguments from II. 1) are still valid. In particular, one has $(\mathcal K(q)f)(\cdot,\cdot)\in BC\left([0,\infty), C_0(l)\right)$ and for $\rho\in l$ inequality {\eqref{est K plus}} is true, that yields the estimate:
\begin{equation}\label{est K plus C0 l}
  \left\|(\mathcal K^+(q)f)(x,\rho)\right\|\leq M  \|q\|_{L_1(0,\infty)} \|f\|_{BC\left([0,\infty), C_0(l)\right)}.
\end{equation}
Furthermore, using the estimates:
\begin{equation}\label{est K plus L2l}
  \left\|(\mathcal K^+(q)f)(x,\cdot)\right\|_{L_2(l)}
  \leq M \int\limits_x^\infty \|q(t)\|\|f(t,\cdot)\|_{L_2(l)} \, dt,
\end{equation}
$$
  \left\|(\mathcal K^+(q)f)(x,\cdot)\right\|_{L_2(l^+(R))}
  \leq M \int\limits_x^\infty \|q(t)\|\|f(t,\cdot)\|_{L_2(l^+(R))} \, dt,
$$
with arbitrary $R>0$ and Lemma 3.1 from {\cite{IgnIntTr}} we obtain that $(\mathcal K^+(q)f)(x,\cdot)\in L_2(l)$ for all $x\in[0,\infty)$, $(\mathcal K^+(q)f)(\cdot,\cdot)\in BC\left([0,\infty), L_2(l)\right)$. Moreover, from {\eqref{est K plus C0 l}} and the estimate
$$
   \left\|(\mathcal K^+(q)f)(x,\cdot)\right\|_{L_2(l)}
  \leq M \|f\|_{BC([0,\infty), L_2(l))} \cdot \|q\|_{L_1(0,\infty)},
$$
it follows that $\mathcal K^+(q)\in \mathcal L\left(BC([0,\infty), \mathcal H(l))\right)$ and is continuous with respect to $q\in\mathcal X_p$.

II. 4) From {\eqref{est K plus L2l}} it follows the estimate:
$$
\left\|((\mathcal K^+(q))^rf)(x,\cdot)\right\|_{L_2(l)}
  \leq \frac{M^r}{r!} \left(\int\limits_x^\infty \|q(t)\|\, dt\right)^r \sup\limits_{\tau\in[x,\infty)}\|f(\tau,\cdot)\|_{L_2(l)},
$$
that yields (together with {\eqref{est Kr plus}}) the estimates similar to
{\eqref{est Kr plus norm}} but with $\mathcal L\left(BC([0,\infty), \mathcal H(l))\right)$ norm in the left hand side. Proceeding as in II. 2), we conclude that
$$\hat F_k\in C\left(\mathcal X_p, BC\left([0,\infty), \mathcal H(l)\right)\right).$$

$\hfil\Box$

\section{Weyl--type solutions}
\label{sec4}

Suppose that $q(\cdot)\in L_1(0,\infty)$, $\rho\in\mathbb C\setminus\Sigma$, $k\in\{1,\dots,n\}$.

\medskip
\textbf{Definition 1.} \textit{Function $y(x)$, $x\in(0,\infty)$ is called {\it $k$-th Weyl-type solution} if it satisfy {\eqref{main sys}} and the following asymptotics hold:
\begin{equation}\label{Weyl def assymp}
  y(x)=O(x^{\mu_k}), x\to 0, \ y(x)=\exp(\rho R_k x)(\mathfrak{f}_k+o(1)), x\to\infty.
\end{equation}
}

Let $\{T_k(\cdot,\rho)\}_{k=1}^n$ and $\{F_k(\cdot,\rho)\}_{k=1}^n$ be the fundamental tensors constructed in the previous section. We define the \textit{characteristic functions} as follows:
$$
\Delta_k(\rho) := |F_{k-1}(x,\rho)\wedge T_k(x,\rho)|, \ k=\overline{2,n}, \ \Delta_1(\rho):=1.
$$

We can repeat the arguments of {\cite{IgnSingRM16}} and obtain the following results similar to {\cite{IgnSingRM16}}, Theorem 4.1 and {\cite{IgnSingRM16}}, Lemma 4.2.

\medskip
\textbf{Proposition 2.} \textit{If $\Delta_k(\rho)\neq 0$ then the $k$-th Weyl type solution exists and is unique. Moreover, it coincides with the (unique for each fixed $x\in (0,\infty)$) solution $\Psi_k(x,\rho)$ of the following SLAE:
$$
F_{k-1}(x,\rho)\wedge \Psi_k(x,\rho) = F_k(x,\rho), \ \Psi_k(x,\rho)\wedge T_k(x,\rho)=0.
$$
}

Our further considerations will be concentrated mostly on the properties of the Weyl type solutions as functions of spectral parameter $\rho$ and potential $q(\cdot)$. In what follows we use the notation $\Psi_k(q,x,\rho)$ for the $k$-th Weyl type solution and the notation $\Psi(q,x,\rho)$ for the matrix $(\Psi_1(q,x,\rho), \dots, \Psi_n(q,x,\rho))$. We also use the notation:
 $\tilde\Psi(q,x,\rho):= \Psi(q,x,\rho) (W(\rho x))^{-1}$, where $W(\xi):=diag(W_1(\xi), \dots, W_n(\xi))$.

As above we fix an arbitrary open sector $\mathcal S\subset\mathbb C\setminus\Sigma$ with a vertex at $\rho=0$.

\medskip
\textbf{Theorem 2.}\textit{
Consider the matrix:
$$
\beta(q,x,\rho):=(\tilde\Psi_0(x,\rho))^{-1}\tilde\Psi(q,x,\rho).
$$
Suppose that $p>2$. Then the following representations hold:
$$
\beta_{jk}(q,x,\rho)=\frac{\delta_{j,k}+d_{jk}(q,x,\rho)}{1+d_k(q,x,\rho)},
$$
where
$\delta_{j,k}$ is the Kronecker delta and $d_{jk}, d_k\in C(\mathcal{X}_p, BC([0,\infty),C_0({\overline{\mathcal{S}}}))$. Moreover, for any ray $l=\{\rho=zt, t\in[0,\infty)\}$, $z\in{\overline{\mathcal{S}}\setminus\{0\}}$  the restrictions $\left.d_{jk}\right|_l, \left.d_k\right|_l$ belong to the space $C(\mathcal{X}_p, BC([0,\infty), \mathcal H(l)))$.
}

\medskip
\textbf{Proof.}
Our definition of the matrix  $\beta(q,x,\rho)$ is equivalent to the following relations:
$$
\tilde\Psi_k(q,x,\rho)=\sum\limits_{j=1}^n \beta_{jk}(q,x,\rho)\tilde\Psi_{0j}(x,\rho), \ k=\overline{1,n}.
$$
Fix an arbitrary $k\in\{1,\dots,n\}$ and some $(q,x,\rho)$. Substituting the relations above into the system:
$$
\tilde F_{k-1} \wedge \tilde \Psi_k = \tilde F_k, \ \tilde\Psi_k \wedge \tilde T_k=0,
$$
where $$\tilde F_k(q,x,\rho):=(\overrightarrow{W}^k(\rho x))^{-1}F_k(q,x,\rho), \ \tilde T_k(q,x,\rho):=(\overleftarrow{W}^k(\rho x))^{-1}T_k(q,x,\rho)$$ we obtain:
$$
\left\{
\begin{array}{l}
  \sum\limits_{j=1}^n \beta_{jk} \tilde F_{k-1}\wedge\tilde\Psi_{0j}=\tilde F_k, \\
  \sum\limits_{j=1}^n \beta_{jk} \tilde\Psi_{0j}\wedge\tilde T_k=0.
\end{array}
\right.
$$
Thus, we arrive at the following SLAE with respect to the values
 $\{\beta_{jk}\}_{j=1}^k$:
\begin{equation}\label{SLAU Psi}
  \sum\limits_{j=1}^n m_{ij} \beta_{jk}=u_i,
\end{equation}
$$
  m_{ij}= \left|\tilde F_{k-1}\wedge\tilde\Psi_{0j}\wedge\tilde\Psi_{0k}
  \wedge\dots\wedge\tilde\Psi_{0,i-1}\wedge\tilde\Psi_{0,i+1}\wedge\dots\wedge\tilde\Psi_{0n}\right|, \ i=\overline{k,n}, $$
$$
  m_{ij}=\left|\tilde\Psi_{01}\wedge\dots\wedge\tilde\Psi_{0,i-1}\wedge\tilde\Psi_{0,i+1}
  \wedge\dots\wedge\tilde\Psi_{0,k-1}\wedge\tilde\Psi_{0,j}\wedge\tilde T_k\right|, \ i=\overline{1,k-1},$$
$$
  u_i= \left|\tilde F_k\wedge\tilde\Psi_{0k}\wedge\dots\wedge\tilde\Psi_{0,i-1}\wedge\tilde\Psi_{0,i+1}\wedge\dots\wedge\tilde\Psi_{0n}\right|, \ i=\overline{k,n}$$
$$  u_i=0, \ i=\overline{1,k-1}.$$
Using Theorem 1 we deduce from the relations above the following representations:
$$
m_{ij}(q,x,\rho)=\left(m^0_{ij}\frac{W_i(\rho x)}{W_j(\rho x)}+\hat m_{ij}(q,x,\rho)\right)\cdot\frac{1}{\overrightarrow{W}^n(\rho x)},$$
$$
  m^0_{ij}= \left|\tilde F^0_{k-1}\wedge\tilde\Psi_{0j}\wedge\tilde\Psi_{0k}
  \wedge\dots\wedge\tilde\Psi_{0,i-1}\wedge\tilde\Psi_{0,i+1}\wedge\dots\wedge\tilde\Psi_{0n}\right|, \ i=\overline{k,n}, $$
$$
  m^0_{ij}=\left|\tilde\Psi_{01}\wedge\dots\wedge\tilde\Psi_{0,i-1}\wedge\tilde\Psi_{0,i+1}
  \wedge\dots\wedge\tilde\Psi_{0,k-1}\wedge\tilde\Psi_{0,j}\wedge\tilde T^0_k\right|, \ i=\overline{1,k-1},$$
$$ u_i(q,x,\rho)=\left(u^0_i\frac{W_i(\rho x)}{W_k(\rho x)}+\hat u_i(q,x,\rho)\right)\cdot\frac{1}{\overrightarrow{W}^n(\rho x)}.
$$
$$
  u^0_i= \left|\tilde F^0_k\wedge\tilde\Psi_{0k}\wedge\dots\wedge\tilde\Psi_{0,i-1}\wedge\tilde\Psi_{0,i+1}\wedge\dots\wedge\tilde\Psi_{0n}\right|, \ i=\overline{k,n}$$
Here the scalar functions $\hat m_{ij}=\hat m_{ij}(q,x,\rho)$, $\hat u_i=\hat u_i(q,x,\rho)$ are such that \\ $\hat m_{ij}, \hat u_i\in C(\mathcal{X}_p, BC([0,\infty),C_0({\overline{\mathcal{S}}}))$ and moreover, for any ray $l=\{\rho=zt, t\in[0,\infty)\}$, $z\in{\overline{\mathcal{S}}}\setminus\{0\}$  the restrictions $\left.\hat m_{ij}\right|_l, \left.\hat u_i\right|_l$ belong to $C(\mathcal{X}_p, BC([0,\infty), \mathcal H(l)))$.

Using the Kramer's rule to solve {\eqref{SLAU Psi}} and taking into account that
$$m^0_{ij}=\delta_{i,j}(-1)^{i-k+1}\Delta^0_{k}, \ i=\overline{1,k-1},$$
$$m^0_{ij}=\delta_{i,j}(-1)^{i-k}|\mathfrak{f}|, \ i=\overline{k,n},$$
$$u_i^0=\delta_{i,k}|\mathfrak{f}|$$
we obtain the desired assertion.

$\hfil\Box$

More detailed results are possible to be obtained if the characteristic function is a priori known to possess only nonzero values on some subset of $\overline{\mathcal S}$.

\medskip
\textbf{Definition 2.} \textit{
Let $L$ be some closed (possibly unbounded) subset of ${\overline{\mathcal{S}}}$. We say that a matrix function $q\in\mathcal X_p$ belongs to $G^p_0(L)$ if
$$\prod\limits_{k=1}^n \Delta_k(\rho)\neq 0$$
for all $\rho\in L$.
}

\medskip
\textbf{Theorem 3.} \textit{
Suppose that $\Psi_k(q,x,\rho)$ is written in the following form:
$$\Psi_k(q,x,\rho)=\Psi_{0k}(x,\rho)+W_k(\rho x)\hat\Psi_k(q,x,\rho).$$
Then for any finite segment $[0,T]$ and any closed subset $L\subset{\overline{\mathcal{S}}}$ it is true that\\ $\hat\Psi_k\in C(G^p_0(L), C([0,T], C_0(L)))$.
}

\medskip
\textbf{Proof.}
1) Using the same arguments as in proof of Theorem 4.4 {\cite{IgnSingRM16}} we can deduce that for any $q\in G^p_0(L)$ the function $\tilde\Psi_k(q,x,\rho)$ admits a continuous extension onto the set $[0,\infty)\times L$. Moreover, for any fixed $(x,\rho)\in[0,\infty)\times L$, $k=\overline{1,n}$ the tensors $\tilde T_k(q,x,\rho)$, $\tilde F_k(q,x,\rho)$ are decomposable.
By virtue of Lemma 7.1 {\cite{BeDT}} the condition $\Delta_k(\rho)\neq 0$, $\rho\in L$ guarantees a unique solvability (for any fixed $(x,\rho)\in [0,\infty)\times L$) of the SLAE:
$$
\tilde F_{k-1} \wedge \tilde \Psi_k = \tilde F_k, \ \tilde\Psi_k \wedge \tilde T_k=0.
$$
Thus, we can conclude that $\tilde\Psi_k\in C(G^p_0(L), C([0,T]\times K))$ for any finite segment $[0,T]$ and any compact set $K\subset L$. This is equivalent to $\tilde\Psi_k\in C(G^p_0(L), C([0,T], C(K)))$. Since the compact set $K$ is \textit{arbitrary} we also have obtained that $\hat\Psi_k(q, \cdot,\cdot)\in C([0,T]\times L)$.

2) It follows from Theorem 2 that for any fixed $q\in\mathcal X_p$ one has $\hat\Psi_k(q,x,\rho)\to 0$ as $\rho\to \infty$ uniformly with respect to $x\in [0,T]$ for any finite $T$. Therefore, we have $\hat\Psi_k(q,\cdot,\cdot)\in C([0,T], C_0(L))$ for any fixed $q\in G^p_0(L)$.

3) Let a finite $T>0$ be arbitrary fixed. Fix an arbitrary  $q_0\in G^p_0(L)$ and consider a sequence $q_m\in G^p_0(L)$ converging to $q_0$ in $\mathcal X_p$ norm.

By virtue of Theorem 2 we have (we use the same notations)
$d_k(q_0,\cdot,\cdot)\in BC([0,\infty), C_0({\overline{\mathcal{S}}}))$.
Therefore, there exists $R>0$ such that $|d_k(q_0,x,\rho)|\geq 1/2$ for all $x\in[0,T]$, $\rho\in L: |\rho|>R$.

Further, from Theorem 2 it follows that:
$$
\sup\limits_{x\in[0,T]}\sup\limits_{\rho\in\overline{\mathcal S},|\rho|>R} |d_k(q_m,x,\rho)-d_k(q_0,x,\rho)|\to 0,$$$$
\sup\limits_{x\in[0,T]}\sup\limits_{\rho\in\overline{\mathcal S},|\rho|>R} |d_{jk}(q_m,x,\rho)-d_{jk}(q_0,x,\rho)|\to 0.
$$
Using this we obtain:
$$
\sup\limits_{x\in[0,T]}\sup\limits_{\rho\in\overline{\mathcal S},|\rho|>R} |\beta_{jk}(q_m,x,\rho)-\beta_{jk}(q_0,x,\rho)|\to 0
$$
as $m\to\infty$. This yields:
$$
\sup\limits_{x\in[0,T]}\sup\limits_{\rho\in\overline{\mathcal S},|\rho|>R}
\|\hat\Psi_{k}(q_m,x,\rho)-\hat\Psi_{k}(q_0,x,\rho)\|\to 0, \ m\to\infty,
$$
while the results obtained in part 1 of this proof imply:
$$
\sup\limits_{x\in[0,T]}\sup\limits_{\rho\in L,|\rho|\leq R}
\|\hat\Psi_{k}(q_m,x,\rho)-\hat\Psi_{k}(q_0,x,\rho)\|\to 0, \ m\to\infty.
$$
Thus, we have:
$$
\lim\limits_{m\to\infty}\|\hat\Psi_{k}(q_m,\cdot,\cdot)-\hat\Psi_{k}(q_0,\cdot,\cdot)\|_{C([0,T], C_0(L))}=0.
$$

$\hfil\Box$

\medskip
\textbf{Theorem 4.} \textit{
Consider the representation:
$$\Psi_k(q,x,\rho)=\Psi_{0k}(x,\rho)+W_k(\rho x)\hat\Psi_k(q,x,\rho),$$
Suppose $p>2$.
Then $\hat\Psi_k\in C(G^p_0(l), C([0,T], \mathcal H(l)))$ for any ray $l=\{\rho=zt, t\in[0,\infty)\}$, $z\in\overline{\mathcal S}\setminus\{0\}$  and any finite segment $[0,T]$.
}

\medskip
\textbf{Proof.}
Fix an arbitrary $T\in(0,\infty)$ and ray $l=\{\rho=zt, t\in[0,\infty)\}$, $z\in\overline{\mathcal S}\setminus\{0\}$. For arbitrary fixed $q\in G^p_0(l)$, $x\in[0,T]$, $\rho\in L$ consider the SLAE:
$$
F_{k-1}\wedge\Psi_k=F_k, \ \Psi_k\wedge T_k=0.
$$
After the substitutions (where $W_j$ assumes $W_j(\rho x)$):
$$
F_{k-1}=F^0_{k-1}+\overrightarrow W^{k-1}\hat F_{k-1}, \ F_k=F^0_k+\overrightarrow W ^k\hat F_k, \
T_k=T^0_k+\overleftarrow W^k \hat T_k,
$$
$$\Psi_k=\Psi_{0k}+W_k\hat\Psi_k$$
the SLAE becomes as follows:
\begin{equation}\label{SLAU 2 Th 2}
\tilde F^0_{k-1}\wedge\hat\Psi_k=f_k, \ \hat\Psi_k\wedge\tilde T^0_k=g_k,
\end{equation}
where:
$$
f_k=\hat F_k -\hat F_{k-1}\wedge\tilde\Psi_{0k}- \hat F_{k-1}\wedge\hat\Psi_k, \
g_k=-\tilde\Psi_{0k}\wedge\hat T_k - \hat\Psi_k \wedge\hat T_k.
$$

From Theorem 1 it follows that the functions $\hat F_{k-1}(\cdot,\cdot,\cdot), \hat F_k(\cdot,\cdot,\cdot), \hat T_k(\cdot,\cdot,\cdot)$ belong to \\ $C(\mathcal X_p, BC([0,\infty), \mathcal H(l)))$.
From Theorem 3 it follows that $\hat\Psi_k(\cdot,\cdot,\cdot)\in C(G^p_0(l), C([0,T], C_0(l)))$. Therefore, the functions $f_k(\cdot,\cdot,\cdot)$ and $g_k(\cdot,\cdot,\cdot)$ belong to $C(G^p_0(l), C([0,T], \mathcal H(l)))$.

We write the function $ \hat\Psi_k(q,x,\rho)$ in the form:
\begin{equation}\label{repr hat Psi k}
  \hat\Psi_k(q,x,\rho)=\sum\limits_{j=1}^n \hat\beta_{jk}(q,x,\rho)\tilde\Psi_{0j}(x,\rho).
\end{equation}
Substituting {\eqref{repr hat Psi k}} into {\eqref{SLAU 2 Th 2}} we arrive after some algebra at the following representations:
$$
\hat\beta_{jk}=$$$$(-1)^{j-k}\left(\Delta_{0k}\overrightarrow{W}^n\right)^{-1}\cdot
\left|\tilde\Psi_{01}\wedge\dots\wedge\tilde\Psi_{0,j-1}\wedge\tilde\Psi_{0,j+1}\wedge\dots\wedge\tilde\Psi_{0,k-1}\wedge g_k\right|, j<k,
$$
$$
\hat\beta_{jk}=$$$$(-1)^{j-k}|\mathfrak{f}|\left(\overrightarrow{W}^n\right)^{-1}\cdot
\left|f_k\wedge\tilde\Psi_{0k}\wedge\dots\wedge\tilde\Psi_{0,j-1}\wedge\tilde\Psi_{0,j+1}\wedge\dots\wedge\tilde\Psi_{0n}\right|, j=\overline{k,n},
$$
that yield
$\hat\beta_{jk}(\cdot,\cdot,\cdot)\in C(G^p_0(l), C([0,T], \mathcal H(l)))$.

$\hfil\Box$

\end{document}